# OPTIMAL INVESTMENT WITH RANDOM ENDOWMENTS IN INCOMPLETE MARKETS

By Julien Hugonnier[1] and Dmitry Kramkov[2]

*HEC Montréal and CIRANO and Carnegie Mellon University*

In this paper, we study the problem of expected utility maximization of an agent who, in addition to an initial capital, receives random endowments at maturity. Contrary to previous studies, we treat as the variables of the optimization problem not only the initial capital but also the number of units of the random endowments. We show that this approach leads to a dual problem, whose solution is always attained in the space of random variables. In particular, this technique does not require the use of finitely additive measures and the related assumption that the endowments are bounded.

**1. Introduction.** A classical problem in financial economics is that of an agent who invests his or her initial capital into a security market so as to maximize the expected utility of terminal wealth. In the context of continuous-time models, this problem was first studied by Merton (1969, 1971). Using dynamic programming arguments, he derived a nonlinear partial differential equation for the value function of this stochastic control problem and obtained closed-form solutions for different specifications of the agent's utility function. The introduction by Harrison and Kreps (1979), Harrison and Pliska (1981) and Ross (1976) of the notion of equivalent martingale measures created the possibility of solving such problems by martingale duality methods. Under the assumption of complete markets, which implies that the family of martingale measures is a singleton, this approach was developed by Cox and Huang (1991, 1998), Karatzas, Lehoczky and Shreve (1987) and Pliska (1986). The essentially more difficult case of incomplete financial markets was considered by He and Pearson (1991a, b),

Received October 2002; revised May 2003.

[1]Research supported by NSF Grant DMS-98-19950.
[2]Research supported by NSF Grant DMS-01-39911.

*AMS 2000 subject classifications.* 90A09, 90A10, 90C26.

*Key words and phrases.* Utility maximization, random endowment, convex duality, incomplete markets, optimal investment, utility-based valuation, contingent claim, European option.

Karatzas, Lehoczky, Shreve and Xu (1991) and more recently by Kramkov and Schachermayer (1999, 2003). In particular, the papers of Kramkov and Schachermayer (1999, 2003) contain minimal conditions on the agent's utility function and the financial market model which imply the key assertions of the theory.

In this paper, we study the expected utility maximization problem of an agent who, in addition to an initial capital, receives random endowments (e.g., the payoffs of contingent claims) at maturity. In complete markets, the agent's random endowments can be perfectly replicated by a controlled portfolio of the traded assets. As a result, the optimal investment problem becomes equivalent to one without random endowments but with an augmented initial capital [see Karatzas and Shreve (1998), Chapter 4]. In incomplete markets, such a transformation becomes impossible in general. In this case, the problem was studied by Cuoco (1997), Cvitanič, Schachermayer and Wang (2001) and Duffie, Fleming, Soner and Zariphopoulou (1997). In particular, Cvitanič, Schachermayer and Wang (2001) consider a general semimartingale model of a financial market and give a characterization of the optimal terminal wealth in terms of the solution to a dual problem. The results of Cvitanić, Schachermayer and Wang (2001) have recently been extended by Karatzas and Žitković (2002) to allow for intertemporal consumption.

The dual problem in Cvitanić, Schachermayer and Wang (2001) is defined over the space $(\mathbf{L}^\infty)^*$ of finitely additive measures. Unfortunately, this formulation requires the stringent assumption that the random endowments are bounded. The novelty of our approach is that we treat as the variables of the optimization problem not only the initial capital as in Cvitanić, Schachermayer and Wang (2001) but also the number of units of random endowments. We show that this extension leads to a dual problem, which does not require the use of finitely additive measures and therefore does not rely on any boundedness assumption.

**2. Main results.** We consider a finite-horizon model of a financial market which consists of $d+1$ assets: a savings account and $d$ stocks. As is common in mathematical finance, we assume that the interest rate is 0; that is, the capital invested into or borrowed from the savings account is constant over time. The price process $S = (S^i)_{1 \leq i \leq d}$ of the stocks is assumed to be a semimartingale on a given filtered probability space $(\Omega, \mathcal{F}, \mathbf{F} = (\mathcal{F}_t)_{0 \leq t \leq T}, \mathbb{P})$, where the filtration $\mathbf{F}$ satisfies the usual conditions of right continuity and completion and $T$ is a finite maturity.

A probability measure $\mathbb{Q}$ is called an *equivalent local martingale measure* if it is equivalent to $\mathbb{P}$ and if $S$ is a local martingale under $\mathbb{Q}$. We denote by $\mathcal{M}$ the family of equivalent local martingale measures and assume that

$$\mathcal{M} \neq \varnothing. \tag{1}$$

This rather mild condition is essentially equivalent to the absence of arbitrage opportunities in the model; see Delbaen and Schachermayer (1994, 1998) for precise statements as well as for further references.

A self-financing portfolio is defined as a pair $(x, H)$, where $x \in \mathbf{R}$ represents the initial capital and $H$ is a predictable $S$-integrable process specifying the number of shares of each stock held in the portfolio. The *wealth* process $X$ of the portfolio evolves in time as the stochastic integral of $H$ with respect to $S$:

$$(2) \qquad X_t \triangleq x + (H \cdot S)_t = x + \int_0^t H_u \, dS_u, \qquad 0 \leq t \leq T.$$

For $x \geq 0$, we denote by $\mathcal{X}(x)$ the set of nonnegative wealth processes whose initial value is equal to $x$, that is,

$$(3) \qquad \mathcal{X}(x) \triangleq \{X \geq 0 : X \text{ satisfies (2) and } X_0 = x\}.$$

Since a process $X \in \mathcal{X}(x)$ is a nonnegative stochastic integral with respect to $S$, it is a local martingale and a global supermartingale under every $\mathbb{Q} \in \mathcal{M}$ [see Ansel and Stricker (1994)].

A nonnegative wealth process in $\mathcal{X}(x)$ is said to be *maximal* if its terminal value cannot be dominated by that of any other process in $\mathcal{X}(x)$. According to Delbaen and Schachermayer (1997) Theorem 2.5, a process is maximal if and only if there exists $\mathbb{Q} \in \mathcal{M}$ under which it is a uniformly integrable martingale. Note that the set of maximal strategies coincides with the set of "good" numéraires in the model [see Delbaen and Schachermayer (1995)].

The nonnegative processes from the set $\mathcal{X}(x)$ constitute the optimization set in the classical problem of optimal investment if the utility function is defined on the positive real line. In the presence of random endowments, one has to extend the domain of the problem and to consider portfolios with possibly negative values. If the endowments are uniformly bounded, as in Cvitanić, Schachermayer and Wang (2001), then the optimization set coincides with the set of *admissible* strategies whose wealth processes are uniformly bounded from below. In the general case of unbounded endowments, the optimization set has to be extended further to the set of *acceptable* strategies.

Following Delbaen and Schachermayer (1997), we say that a wealth process $X$ is *acceptable* if it admits a representation of the form $X = X' - X''$, where $X'$ is a nonnegative wealth process and $X''$ is a maximal process. Note that if the maximal process $1 + X''$ is chosen as a numéraire, then the discounted process $X/(1 + X'')$ is uniformly bounded from below and hence is admissible under this numéraire. The definition of acceptable strategies is therefore very natural: they represent the numéraire-invariant version of admissible strategies.

We also consider an economic agent whose preferences over terminal consumption bundles are represented by a utility function $U:(0,\infty) \to \mathbf{R}$. The function $U$ is assumed to be strictly concave, strictly increasing and continuously differentiable and to satisfy the Inada conditions:

$$(4) \qquad U'(0) \triangleq \lim_{x \to 0} U'(x) = \infty, \qquad U'(\infty) \triangleq \lim_{x \to \infty} U'(x) = 0.$$

The convex conjugate of the agent's utility function is defined as follows:

$$V(y) \triangleq \sup_{x>0}\{U(x) - xy\}, \qquad y > 0.$$

It is well known that, under the Inada conditions (4), the conjugate of $U$ is a continuously differentiable, strictly decreasing and strictly convex function satisfying $V'(0) = -\infty$, $V'(\infty) = 0$ and $V(0) = U(\infty)$, $V(\infty) = U(0)$ as well as the following bidual relation:

$$U(x) = \inf_{y>0}\{V(y) + xy\}, \qquad x > 0.$$

Assume now that the portfolio of the agent at time 0 consists of an initial capital $x$ as well as of the quantities $q = (q_i)_{1 \leq i \leq N}$ of *nontraded* European contingent claims with maturity $T$ and $\mathcal{F}_T$-measurable payment functions $f = (f_i)_{1 \leq i \leq N}$. We denote by

$$\langle q, f \rangle \triangleq \sum_{i=1}^{N} q_i f_i$$

the payoff of this portfolio of the contingent claims and by

$$\mathcal{X}(x,q) \triangleq \{X : X \text{ is acceptable, } X_0 = x \text{ and } X_T + \langle q, f \rangle \geq 0\}$$

the set of acceptable processes with initial capital $x$ whose terminal value dominates the random payoff $-\langle q, f \rangle$. Note that in the case when $x > 0$ and $q = 0$ this set coincides with the set of nonnegative wealth processes with initial capital $x$. In other words, we have $\mathcal{X}(x,0) = \mathcal{X}(x)$ for all $x > 0$.

The family $\mathcal{X}(x,q)$ may very well be empty for some vector $(x,q)$. From now on, we shall restrict ourselves to the set $\mathcal{K}$, which is defined as the *interior* of the set of points $(x,q)$, where $\mathcal{X}(x,q)$ is not empty:

$$\mathcal{K} \triangleq \text{int}\{(x,q) \in \mathbf{R}^{N+1} : \mathcal{X}(x,q) \neq \varnothing\}.$$

As is easily seen, this set is a convex cone in $\mathbf{R}^{N+1}$. Hereafter we shall assume that it contains any point $(x,q)$ such that $x > 0$ and $q = 0$, that is

$$(5) \qquad (x,0) \in \mathcal{K}, \qquad x > 0.$$

The condition above means, in particular, that the optimization problem with random endowments contains the classical problem of optimal investment as a special case. The following lemma provides a list of equivalent assertions to (5).

LEMMA 1. *Assume that* (1) *holds true. Then the following conditions are equivalent:*

(i) *The set $\mathcal{K}$ satisfies* (5).
(ii) *For any $q \in \mathbf{R}^N$, there exists $x \in \mathbf{R}$ such that the set $\mathcal{X}(x,q)$ is not empty.*
(iii) *There is a nonnegative wealth process $X$ such that $X_T \geq \sum_{i=1}^{d} |f_i|$.*
(iv) *The payment functions $f = (f_i)_{1 \leq i \leq N}$ satisfy the integrability conditions:*

(6) $$\sup_{Q \in \mathcal{M}} \mathbb{E}_Q[|f_i|] < \infty \qquad \text{for all } 1 \leq i \leq N.$$

PROOF. The equivalences among (i)–(iii) are straightforward, while the equivalence between (iii) and (iv) follows from the general duality relationships between terminal capitals and martingale measures [see, e.g., Delbaen and Schachermayer (1994), Theorem 5.7]. □

The quantity $q_i$ of each contingent claim in the agent's portfolio being held constant up to maturity, the vector $q$ represents the *illiquid* part of the portfolio. On the contrary, the initial wealth $x$ can be freely invested into the stocks and the savings account according to some dynamic strategy. Given $x$ and $q$, the goal of the agent is to maximize the expected utility of his or her terminal wealth. This leads to the following optimization problem:

(7) $$u(x,q) \triangleq \sup_{X \in \mathcal{X}(x,q)} \mathbb{E}[U(X_T + \langle q, f \rangle)], \qquad (x,q) \in \mathcal{K}.$$

REMARK 1. Contrary to Cvitanić, Schachermayer and Wang (2001) and Karatzas and Žitković (2002), we consider $u$ to be the function of *both* $x$ and $q$. As we shall see later, this approach will permit us to avoid the use of finitely additive measures in the formulation of the dual problem and to overcome the related boundedness assumption on $f$.

Note that such a definition of the value function is also useful for the study of utility-based valuation of contingent claims [see Davis (1997), Frittelli (2000) and Hodges and Neuberger (1989)]. For example, a *certainty equivalence price* for the contingent claims $f$ given a portfolio $(x,q)$ is defined as a vector $\widehat{p}(x,q) \in \mathbf{R}^N$ such that

$$w(x + \langle q, \widehat{p}(x,q) \rangle) = u(x,q),$$

where

$$w(x) \triangleq u(x,0) = \sup_{X \in \mathcal{X}(x)} \mathbb{E}[U(X_T)]$$

is the value function of the problem of optimal investment without random endowments. Further, a *utility-based price* for $f$ given $(x,q)$ is defined as a vector $\widetilde{p}(x,q) \in \mathbf{R}^N$ such that the agent's holdings $q$ in the claims are optimal in the model where the claims can be traded at time 0 at price $\widetilde{p}(x,q)$. In other words,

$$u(x,q) \geq u(x',q') \quad \text{if } (x',q') \in \mathcal{K} \quad \text{and} \quad x + \langle q, \widetilde{p}(x,q) \rangle = x' + \langle q', \widetilde{p}(x,q) \rangle.$$

Using standard arguments from the theory of convex functions, we deduce that a vector $\widetilde{p}(x,q) \in \mathbf{R}^N$ is a utility-based price for $f$ given $(x,q)$ if and only if it has the representation:

$$\widetilde{p}(x,q) = \frac{r}{y} \quad \text{for some } (y,r) \in \partial u(x,q),$$

where $\partial u(x,q)$ is the subdifferential of $u$ at the point $(x,q)$. In particular, we deduce that the utility-based price is unique at $(x,q)$ if and only if the value function $u$ is differentiable at that point.

We start our construction of a dual problem to (7) by introducing the set $\mathcal{L}$, which is the *relative interior* of the polar cone of $-\mathcal{K}$:

(8) $\qquad \mathcal{L} \triangleq \text{ri}\{(y,r) \in \mathbf{R}^{N+1} : xy + \langle q, r \rangle \geq 0 \text{ for all } (x,q) \in \mathcal{K}\}.$

The set $\mathcal{L}$ will define the domain of the dual problem at time 0. It can be shown that the intersection of $\mathcal{L}$ with the hyperplane $y \equiv 1$, that is,

(9) $\qquad\qquad\qquad \mathcal{P} = \{p \in \mathbf{R}^N : (1,p) \in \mathcal{L}\},$

defines the set of arbitrage-free prices of the contingent claims $f$.

Following Kramkov and Schachermayer (1999), we denote by $\mathcal{Y}(y)$ the family of nonnegative semimartingales $Y$ with initial value $y$ and such that for any positive wealth process $X$ the product $XY$ is a nonnegative supermartingale, that is,

(10) $\mathcal{Y}(y) \triangleq \{Y \geq 0 : Y_0 = y, XY \text{ is a supermartingale for all } X \in \mathcal{X}(1)\}.$

In particular, as $\mathcal{X}(1)$ contains the constant process 1, the elements of $\mathcal{Y}(y)$ are nonnegative supermartingales. Note also that the set $\mathcal{Y}(1)$ contains the density processes of all $\mathbb{Q} \in \mathcal{M}$.

Given an arbitrary vector $(y,r)$ in $\mathcal{L}$, we denote by $\mathcal{Y}(y,r)$ the set of nonnegative supermartingales $Y \in \mathcal{Y}(y)$ such that the inequality

(11) $\qquad\qquad\qquad \mathbb{E}[Y_T(X_T + \langle q, f \rangle)] \leq xy + \langle q, r \rangle$

holds true for all $(x,q) \in \mathcal{K}$ and $X \in \mathcal{X}(x,q)$. Using this notation, we now define the dual optimization problem to (7) as follows:

(12) $\qquad\qquad v(y,r) \triangleq \inf_{Y \in \mathcal{Y}(y,r)} \mathbb{E}[V(Y_T)], \qquad (y,r) \in \mathcal{L}.$

The following theorems constitute our main results.

THEOREM 1. *Assume that conditions* (1), (4) *and* (5) *hold true and*

(13) $$u(x,q) < \infty \quad \text{for some } (x,q) \in \mathcal{K}.$$

*Then we have:*

(i) *The function $u$ is finitely valued on $\mathcal{K}$ and for any $(y,r) \in \mathcal{L}$ there exists a constant $c = c(y,r) > 0$ such that $v(cy, cr)$ is finite. The value functions $u$ and $v$ are conjugate:*

(14)
$$u(x,q) = \inf_{(y,r)\in\mathcal{L}} \{v(y,r) + xy + \langle q, r\rangle\}, \qquad (x,q) \in \mathcal{K},$$
$$v(y,r) = \sup_{(x,q)\in\mathcal{K}} \{u(x,q) - xy - \langle q, r\rangle\}, \qquad (y,r) \in \mathcal{L}.$$

(ii) *The solution $\widehat{Y}(y,r)$ to* (12) *exists and is unique for all $(y,r) \in \mathcal{L}$ such that $v(y,r) < \infty$.*

REMARK 2. It might be convenient to verify the validity of (13) at a point $(x,q)$, where $x > 0$ and $q = 0$ as in this case

$$u(x,0) = w(x) \triangleq \sup_{X\in\mathcal{X}(x)} \mathbb{E}[U(X_T)]$$

is the value function of the problem of optimal investment without endowments.

THEOREM 2. *Assume that conditions* (1), (4) *and* (5) *hold true and*

(15) $$v(y,r) < \infty \quad \text{for all } (y,r) \in \mathcal{L}.$$

*Then, in addition to the assertions of Theorem* 1, *we have:*

(i) *The subdifferential of $u$ maps $\mathcal{K}$ into $\mathcal{L}$, that is,*

(16) $$\partial u(x,q) \subset \mathcal{L}, \qquad (x,q) \in \mathcal{K}.$$

(ii) *The solution $\widehat{X}(x,q)$ to* (7) *exists and is unique for any $(x,q) \in \mathcal{K}$. In addition, if $(y,r) \in \partial u(x,q)$, then the terminal values of the corresponding solutions are related by*

(17) $$\widehat{Y}_T(y,r) = U'(\widehat{X}_T(x,q) + \langle q, f\rangle),$$

(18) $$\mathbb{E}[\widehat{Y}_T(y,r)(\widehat{X}_T(x,q) + \langle q, f\rangle)] = xy + \langle q, r\rangle.$$

REMARK 3. The relationship (16) has the economical interpretation that the utility-based prices defined in Remark 1 are arbitrage-free prices for the contingent claims $f$.

It can be shown that contrary to (16) the subdifferential of $-v$ at $(y,r) \in \mathcal{L}$ is not necessarily contained in $\mathcal{K}$. However, it is always contained in the closure of this set.

The proofs of Theorems 1 and 2 will be given in Section 3. The validity of condition (15) might be difficult to verify directly. The following lemma provides an equivalent condition in terms of the function

$$\widetilde{w}(y) \triangleq \inf_{Y \in \mathcal{Y}(y)} \mathbb{E}[V(Y_T)], \qquad y > 0,$$

which is the value function of the dual problem in the classical problem of optimal investment without random endowments [see Kramkov and Schermayer (1999, 2003)]. Recall that we denote by $\mathcal{P}$ the intersection of $\mathcal{L}$ with the hyperplane $y \equiv 1$; see (9).

LEMMA 2. *Assume that conditions* (1), (4) *and* (5) *hold true. Then* (15) *holds if and only if*

$$\widetilde{w}(y) < \infty, \qquad y > 0.$$

*Moreover, in this case,*

(19) $$\widetilde{w}(y) = \inf_{p \in \mathcal{P}} v(y, yp),$$

*where the lower bound is attained.*

PROOF. Define the function

$$\widehat{w}(y) = \inf_{p \in \mathcal{P}} v(y, yp), \qquad y > 0.$$

The function $\widehat{w}$ is clearly convex and by Theorem 1 it satisfies

(20) $$w(x) \triangleq u(x, 0) = \inf_{y > 0} \{\widehat{w}(y) - xy\}, \qquad x > 0.$$

According to Theorem 2.1 in Kramkov and Schachermayer (1999), the function $\widetilde{w}$ satisfies a similar equality:

$$w(x) = \inf_{y > 0} \{\widetilde{w}(y) - xy\}, \qquad x > 0.$$

It follows that the functions $\widehat{w}$ and $\widetilde{w}$ are equal to each other, that is, (19) holds true.

If the function $v$ is finitely valued, then so clearly is $\widehat{w} = \widetilde{w}$. Conversely, if $\widehat{w}$ is finitely valued, then the closure of the set $A$ defined by

$$A = \{(y, r) \in \mathcal{L} : v(y, r) < \infty\}$$

contains the origin. Moreover, the set $A$ is convex, $\lambda A \subset A$ if $\lambda \geq 1$ and, by Theorem 1(i),

$$\mathcal{L} = \bigcup_{\lambda > 0} \lambda A.$$

This readily implies that $A = \mathcal{L}$, that is, that the function $v$ is finitely valued on $\mathcal{L}$.

Finally, fix $y > 0$ and let $x = -\widehat{w}'(y)$. From (20), we deduce that

$$w'(x) = y, \qquad w(x) + xy = \widehat{w}(y).$$

From Theorem 2(ii), we deduce the existence of $\widehat{p} \in \mathcal{P}$ such that $(y, y\widehat{p}) \in \partial u(x, 0)$. It follows that

$$v(y, y\widehat{p}) = u(x, 0) + xy = w(x) + xy = \widehat{w}(y). \qquad \square$$

Following Kramkov and Schachermayer (1999), we also formulate a convenient sufficient condition for the validity of the assertions of Theorem 2. We recall that the asymptotic elasticity of the utility function $U$ is defined to be

$$\mathrm{AE}(U) \triangleq \limsup_{x \to \infty} \frac{xU'(x)}{U(x)}.$$

COROLLARY 1. *In addition to the conditions of Theorem 1, assume that*

$$\mathrm{AE}(U) < 1. \tag{21}$$

*Then* (15) *and all the assertions of Theorem 2 hold true.*

PROOF. The condition $\mathrm{AE}(U) < 1$ is equivalent to the following property of $V$ [see Lemma 6.3 in Kramkov and Schachermayer (1999)]: for any constant $c > 0$, there are positive constants $c_1$ and $c_2$ such that

$$V(y/c) \le c_1 V(y) + c_2, \qquad y > 0.$$

Let us fix $(y, r) \in \mathcal{L}$. Since, for any $c > 0$,

$$\mathcal{Y}(y, r) = \mathcal{Y}(cy, cr)/c,$$

we deduce that $v(y, r) < \infty$ if there is a constant $c = c(y, r) > 0$ such that $v(cy, cr) < \infty$. However, the existence of such a constant has been established in Theorem 1. $\square$

REMARK 4. Note that (21) is the *minimal* condition on the utility function only, which implies the finiteness of the dual function $v$ and the assertions of Theorem 2. See the counterexamples in Kramkov and Schachermayer (1999) for the case where there are no random endowments.

Although the functions $U$ and $-V$ are strictly concave and continuously differentiable, the value functions $u$ and $-v$ do not necessarily inherit any of these properties. Note that, due to the conjugacy relations (14), the continuous differentiability of one of these two functions is closely related to the strict concavity of the other. The following easy lemma provides a set of necessary and sufficient conditions for $u$ to be strictly concave and dually for $v$ to be continuously differentiable.

Recall that a random variable $g$ is *replicable* if there is an acceptable process $X$ such that $-X$ is also acceptable and $X_T = g$. Provided that it exists, such a process $X$ is unique and is called the *replication* process for $g$.

LEMMA 3. *Assume that the conditions of Theorem 2 hold true. Then the following assertions are equivalent:*

(i) *The function $u$ is strictly concave on $\mathcal{K}$.*
(ii) *The function $v$ is continuously differentiable on $\mathcal{L}$.*
(iii) *For any $q \in \mathbf{R}^N$ such that $q \neq 0$, the random variable $\langle q, f \rangle$ is not replicable.*

PROOF. The equivalence between (i) and (ii) is a well-known consequence of (14) and (16); see, for example, Theorems 4.1.1 and 4.1.2. in Hiriart-Urruty and Lemaréchal (2001). Further, since $U$ is a strictly concave function, the value function $u$ is strictly concave if and only if for any two distinct points $(x_i, q_i) \in \mathcal{K}$, $i = 1, 2$, the terminal capitals of the corresponding optimal strategies are different:

$$\mathbb{P}[\widehat{X}_T(x_1, q_1) + \langle q_1, f_1 \rangle \neq \widehat{X}_T(x_2, q_2) + \langle q_2, f_2 \rangle] > 0.$$

As is easily seen, this last property is equivalent to condition (iii). □

REMARK 5. Lemma 7 provides two additional conditions, which are equivalent to the assertions of Lemma 3.

The situation with the continuous differentiability of $u$ and the strict convexity of $v$ is more complicated. As mentioned in Remark 1, this property is equivalent to the uniqueness of the utility-based price. This question is studied in our joint paper with Hugonnier, Kramkov and Schachermayer (2003), where we show that, in general, the utility-based price may not be unique and hence the function $u$ may not be differentiable.

**3. Proofs of the main theorems.** The proofs of Theorems 1 and 2 are based on Proposition 1. For arbitrary vectors $(x, q) \in \mathcal{K}$ and $(y, r) \in \mathcal{L}$, denote

(22) $\quad \mathcal{C}(x, q) \triangleq \{g \in \mathbf{L}_+^0 : g \leq X_T + \langle q, f \rangle \text{ for some } X \in \mathcal{X}(x, q)\},$

(23) $\quad \mathcal{D}(y, r) \triangleq \{h \in \mathbf{L}_+^0 : h \leq Y_T \text{ for some } Y \in \mathcal{Y}(y, r)\}.$

With this notation, the value functions $u$ and $v$ defined in (7) and (12) take the form:

$$u(x, q) = \sup_{g \in \mathcal{C}(x,q)} \mathbb{E}[U(g)], \tag{24}$$

$$v(y, r) = \inf_{h \in \mathcal{D}(y,r)} \mathbb{E}[V(h)]. \tag{25}$$

PROPOSITION 1. *Assume that conditions* (1) *and* (5) *hold true. Then the families* $(\mathcal{C}(x,q))_{(x,q) \in \mathcal{K}}$ *and* $(\mathcal{D}(y,r))_{(y,r) \in \mathcal{L}}$ *defined in* (22) *and* (23) *have the following properties:*

(i) *For any* $(x, q) \in \mathcal{K}$, *the set* $\mathcal{C}(x, q)$ *contains a strictly positive constant. A nonnegative function* $g$ *belongs to* $\mathcal{C}(x, q)$ *if and only if*

$$\mathbb{E}[gh] \leq xy + \langle q, r \rangle \quad \text{for all } (y, r) \in \mathcal{L} \text{ and } h \in \mathcal{D}(y, r). \tag{26}$$

(ii) *For any* $(y, r) \in \mathcal{L}$, *the set* $\mathcal{D}(y, r)$ *contains a strictly positive random variable. A nonnegative function* $h$ *belongs to* $\mathcal{D}(y, r)$ *if and only if*

$$\mathbb{E}[gh] \leq xy + \langle q, r \rangle \quad \text{for all } (x, q) \in \mathcal{K} \text{ and } g \in \mathcal{C}(x, q). \tag{27}$$

The proof of Proposition 1 will be broken into several lemmas. According to Lemma 1(iii), there is a maximal positive wealth process $X'$ whose terminal capital dominates the sum of the absolute values of the endowments:

$$X'_T \geq \sum_{i=1}^{N} |f_i|.$$

Denote by $\mathcal{M}'$ the set of equivalent local martingale measures $\mathbb{Q}$ such that $X'$ is a uniformly integrable martingale under $\mathbb{Q}$. We have that $\mathcal{M}'$ is a nonempty, convex subset of $\mathcal{M}$ which is dense in $\mathcal{M}$ with respect to the variation norm [see Delbaen and Schachermayer (1997), Theorem 5.2].

LEMMA 4. *Assume that the conditions of Proposition* 1 *hold true. Let* $(x, q) \in \mathcal{K}$, $X \in \mathcal{X}(x, q)$ *and* $\mathbb{Q} \in \mathcal{M}'$. *Then* $X$ *is a supermartingale under* $\mathbb{Q}$.

PROOF. Denote

$$c \triangleq \max_{1 \leq i \leq N} |q_i|, \qquad Z \triangleq X + cX'.$$

Since $X$ and $X'$ are acceptable processes, $Z$ is an acceptable process. In addition, the terminal value of $Z$ is nonnegative, because

$$Z_T = X_T + cX'_T \geq X_T + \langle q, f \rangle \geq 0.$$

It follows that $Z$ is a nonnegative wealth process and hence is a supermartingale under $\mathbb{Q}$. Since $X'$ is a uniformly integrable martingale under $\mathbb{Q}$, we deduce that $X = Z - cX'$ is a supermartingale under $\mathbb{Q}$. $\square$

The following lemma is a variant of duality relations between the terminal values of wealth processes and martingale measures.

LEMMA 5. *Assume that the conditions of Proposition* 1 *hold true. Let $g$ be a random variable such that $g \geq -cX'_T$, where $c \geq 0$ is a positive constant, and denote*

$$\alpha(g) \triangleq \sup_{\mathbb{Q} \in \mathcal{M}'} \mathbb{E}_{\mathbb{Q}}[g] < \infty. \tag{28}$$

*Then there is an acceptable process $X$ such that $X_0 = \alpha(g)$ and $X_T \geq g$.*

PROOF. Denote

$$h \triangleq g + cX'_T, \qquad \alpha(h) \triangleq \sup_{\mathbb{Q} \in \mathcal{M}'} \mathbb{E}_{\mathbb{Q}}[h].$$

We have $h \geq 0$, $\alpha(h) = \alpha(g) + cX'_0$ and

$$\alpha(h) = \sup_{\mathbb{Q} \in \mathcal{M}} \mathbb{E}_{\mathbb{Q}}[h],$$

because $\mathcal{M}'$ is dense in $\mathcal{M}$ with respect to the variation norm. From the duality relations between the terminal values of wealth processes and martingale measures [see Delbaen and Schachermayer (1998), Theorem 5.12], we deduce the existence of a nonnegative wealth process $Z$ such that $Z_0 = \alpha(h)$ and $Z_T \geq h$. Clearly, $X \triangleq Z - cX'$ satisfies the assertions of the lemma. $\square$

LEMMA 6. *Assume that the conditions of Proposition* 1 *hold true. Then*

$$\operatorname{cl} \mathcal{K} = \{(x,q) \in \mathbf{R}^{N+1} : \mathcal{X}(x,q) \neq \varnothing\},$$

*where $\operatorname{cl} \mathcal{K}$ denotes the closure of the set $\mathcal{K}$ in $\mathbf{R}^{N+1}$.*

PROOF. Let $(x,q) \in \operatorname{cl} \mathcal{K}$ and let $(x^n, q^n)_{n \geq 1}$ be a sequence in $\mathcal{K}$ that converges to $(x,q)$. The assertion of the lemma will follow if we show that $\mathcal{X}(x,q) \neq \varnothing$. Fix $X^n \in \mathcal{X}(x^n, q^n)$, $n \geq 1$, and denote

$$c \triangleq \sup_{n \geq 1} \max_{1 \leq i \leq N} |q_i^n|.$$

Let $\tau$ be a dense countable subset of $\mathbf{R}_+$. The processes $Z^n \triangleq X^n + cX'$, $n \geq 1$, are nonnegative supermartingales under any $\mathbb{Q} \in \mathcal{M}$, and passing if necessary to convex combinations [see Lemma 5.2 in Föllmer and Kramkov

(1997)], we can assume that they Fatou converge on $\tau$ to a process $Z$, that is,

$$Z_t = \limsup_{s \downarrow t, s \in \tau} \limsup_{n \to \infty} Z_s^n = \liminf_{s \downarrow t, s \in \tau} \liminf_{n \to \infty} Z_s^n.$$

Clearly,

$$\begin{aligned} Z_T - cX_T' + \langle q, f \rangle &= \lim_{n \to \infty} (Z_T^n - cX_T' + \langle q^n, f \rangle) \\ &= \lim_{n \to \infty} (X_T^n + \langle q^n, f \rangle) \geq 0. \end{aligned}$$

The process $Z$ being a nonnegative supermartingale under any $\mathbb{Q} \in \mathcal{M}$, it admits an optional decomposition: $Z = Z' - C$ where $Z'$ is a nonnegative wealth process and $C$ is an increasing process with initial value 0 [see Kramkov (1996)]. Moreover, we deduce from Fatou's lemma that

$$Z_0 \leq \liminf_{n \to \infty} Z_0^n = x + cX_0'.$$

It follows that the process

$$X = Z' - cX' + (x + cX_0' - Z_0)$$

belongs to $\mathcal{X}(x, q)$ and hence the set $\mathcal{X}(x, q)$ is not empty. $\square$

LEMMA 7. *Assume that the conditions of Proposition* 1 *hold true. Then the following statements are equivalent:*

(i) *The set $\mathcal{L}$ is open in $\mathbf{R}^{N+1}$.*

(ii) *For any $q \in \mathbf{R}^N$ such that $q \neq 0$, the random variable $\langle q, f \rangle$ is not replicable.*

(iii) *For any nonzero vector $(x, q) \in \mathrm{cl}\,\mathcal{K}$, there is $X \in \mathcal{X}(x, q)$ such that $\mathbb{P}[X_T + \langle q, f \rangle > 0] > 0$.*

PROOF. By the properties of polars of convex sets [see Rockafellar (1970), Corollary 14.6.1], the first item of the lemma is equivalent to the fact that the set $\mathrm{cl}\,\mathcal{K}$ does not contain any lines passing through the origin, that is,

(29) $\quad (x, q) \neq 0, \quad (x, q) \in \mathrm{cl}\,\mathcal{K} \Rightarrow (-x, -q) \notin \mathrm{cl}\,\mathcal{K},$

where $\mathrm{cl}\,\mathcal{K}$ denotes the closure of $\mathcal{K}$. Now the equivalence of all the statements of the lemma is a rather straightforward consequence of Lemma 6. $\square$

REMARK 6. In some of the proofs below, we shall assume that the equivalent assertions of Lemma 7 hold true. This is without loss of generality. Indeed, we can always arrive at this case by replacing the original family $f$ with a subset $f'$ chosen in such a way that:

1. Any element of $f$ is a linear combination of elements of $f'$ and terminal capitals of acceptable strategies.
2. Any nonzero linear combination of elements of $f'$ is not replicable.

Note that the set $f'$ is empty if and only if all the elements of $f$ are replicable.

For a vector $p \in \mathbf{R}^N$, we denote by $\mathcal{M}'(p)$ the subset of $\mathcal{M}'$ that consists of measures $\mathbb{Q} \in \mathcal{M}'$ such that $\mathbb{E}_\mathbb{Q}[f] = p$. Recall that the set $\mathcal{P}$ denotes the intersection of $\mathcal{L}$ with the hyperplane $y \equiv 1$; see (9).

LEMMA 8. *Assume that the conditions of Proposition* 1 *hold true and let* $p \in \mathbf{R}^N$. *Then the set* $\mathcal{M}'(p)$ *is not empty if and only if* $p \in \mathcal{P}$. *In particular,*

$$\bigcup_{p \in \mathcal{P}} \mathcal{M}'(p) = \mathcal{M}'. \tag{30}$$

PROOF. Define the set

$$\mathcal{P}' \triangleq \{p \in \mathbf{R}^N : \mathcal{M}'(p) \neq \varnothing\}.$$

We have to prove that $\mathcal{P} = \mathcal{P}'$. Since both of these sets are convex and $\mathcal{P}$ is relatively open, we have that $\mathcal{P}$ is contained in $\mathcal{P}'$ if and only if

$$\sup_{p \in \mathcal{P}} \langle q, p \rangle \leq \sup_{p \in \mathcal{P}'} \langle q, p \rangle \qquad \text{for all } q \in \mathbf{R}^N. \tag{31}$$

Fix $q \in \mathbf{R}^N$ and denote

$$\beta(q) \triangleq \sup_{p \in \mathcal{P}'} \langle q, p \rangle = \sup_{\mathbb{Q} \in \mathcal{M}'} \mathbb{E}_\mathbb{Q}[\langle q, f \rangle].$$

From Lemma 5, we deduce the existence of an acceptable process $X$ such that $X_0 = \beta(q)$ and $X_T \geq \langle q, f \rangle$. By Lemma 6, this implies that $(\beta(q), -q) \in \text{cl}\,\mathcal{K}$ and hence that

$$\beta(q) - \langle q, p \rangle \geq 0 \qquad \text{for any } p \in \mathcal{P},$$

which proves (31).

For the proof of the inverse inclusion, it is convenient to assume that the assertions of Lemma 7 are valid. As was noted in Remark 6, this is without any loss of generality. Let $p \in \mathcal{P}'$, $(x, q) \in \text{cl}\,\mathcal{K}$, $\mathbb{Q} \in \mathcal{M}'(p)$ and $X \in \mathcal{X}(x, q)$ be as in assertion (iii) of Lemma 7, that is, such that

$$\mathbb{P}[X_T + \langle q, f \rangle > 0] > 0.$$

Taking into account the supermartingale property of $X$ under $\mathbb{Q}$ established in Lemma 4, we deduce that

$$0 < \mathbb{E}_\mathbb{Q}[X_T + \langle q, f \rangle] \leq x + \langle q, p \rangle.$$

As $(x, q)$ is an arbitrary element of $\text{cl}\,\mathcal{K}$, this implies that $p \in \mathcal{P}$. □

LEMMA 9. *Assume that the conditions of Proposition* 1 *hold true and let $p \in \mathcal{P}$. Then the density process of any $\mathbb{Q} \in \mathcal{M}'(p)$ belongs to $\mathcal{Y}(1,p)$.*

PROOF. The result follows from the supermartingale characterization of wealth processes provided by Lemma 4. □

LEMMA 10. *Assume that the conditions of Proposition* 1 *hold true. For any $(x,q) \in \mathcal{K}$, a nonnegative function $g$ belongs to $\mathcal{C}(x,q)$ if and only if*

(32) $\qquad \mathbb{E}_{\mathbb{Q}}[g] \leq x + \langle q,p \rangle \qquad \text{for all } p \in \mathcal{P} \quad \text{and} \quad \mathbb{Q} \in \mathcal{M}'(p).$

PROOF. If $g \in \mathcal{C}(x,q)$, then the validity of (32) follows from Lemma 9. For the converse, assume that $g$ is a nonnegative random variable such that (32) holds true and denote

$$h \triangleq g - \langle q, f \rangle, \qquad c \triangleq \max_{1 \leq i \leq N} |q_i|.$$

We have $h \geq -cX'_T$ and

$$\alpha(h) \triangleq \sup_{\mathbb{Q} \in \mathcal{M}'} \mathbb{E}_{\mathbb{Q}}[h] = \sup_{p \in \mathcal{P}} \sup_{\mathbb{Q} \in \mathcal{M}'(p)} \mathbb{E}_{\mathbb{Q}}[h]$$
$$= \sup_{p \in \mathcal{P}} \sup_{\mathbb{Q} \in \mathcal{M}'(p)} (\mathbb{E}_{\mathbb{Q}}[g] - \langle q, p \rangle) \leq x,$$

where in the second equality we used (30). Lemma 5 implies the existence of an acceptable process $X$ such that $X_0 = \alpha(h)$ and $X_T \geq h$. It follows that

$$X_T + \langle q, f \rangle \geq g \geq 0.$$

Therefore, $X$ belongs to $\mathcal{X}(\alpha(h), q)$ and $g$ is an element of $\mathcal{C}(x,q)$. □

We are now able to complete the proof of Proposition 1.

PROOF OF PROPOSITION 1. We start with assertion (i). Let $(x,q) \in \mathcal{K}$. Since $\mathcal{K}$ is an open set, there is $\delta > 0$ such that $(x - \delta, q) \in \mathcal{K}$. If $X$ lies in $\mathcal{X}(x - \delta, q)$, then $Z = X + \delta$ belongs to $\mathcal{X}(x,q)$ and

$$Z_T + \langle q, f \rangle \leq \delta.$$

It follows that $\delta \in \mathcal{C}(x,q)$. If $g \in \mathcal{C}(x,q)$, then (26) follows from the construction of the sets $\mathcal{D}(y,r)$, $(y,r) \in \mathcal{L}$. Conversely, let us assume that $g$ is a nonnegative random variable such that (26) holds true. By Lemma 9, the density processes of $\mathbb{Q} \in \mathcal{M}'(p)$ belong to $\mathcal{Y}(1,p)$ for all $p \in \mathcal{P}$ and hence $g$ satisfies (32). Lemma 10 implies now that $g \in \mathcal{C}(x,q)$ and the proof of assertion (i) is complete.

We now turn to assertion (ii). As

$$c\mathcal{D}(y,r) = \mathcal{D}(cy, cr) \qquad \text{for all } c > 0, (y,r) \in \mathcal{L},$$

we can restrict ourselves to the case where $(y,r) = (1,p)$ for some $p \in \mathcal{P}$. From Lemma 8, we deduce the existence of $\mathbb{Q} \in \mathcal{M}'(p)$. By Lemma 9, the Radon–Nikodym derivative $d\mathbb{Q}/d\mathbb{P}$ belongs to $\mathcal{D}(1,p)$. Since the measures $\mathbb{Q}$ and $\mathbb{P}$ are equivalent, we have $\mathbb{P}(d\mathbb{Q}/d\mathbb{P} > 0) = 1$.

If $h \in \mathcal{D}(1,p)$, then (27) follows from the construction of the set $\mathcal{Y}(1,p)$. Conversely, let us assume that $h$ is a nonnegative random variable such that (27) holds true. Then, in particular,

$$\mathbb{E}[gh] \leq 1 \qquad \text{for all } g \in \mathcal{C}(1,0).$$

Proposition 3.1 in Kramkov and Schachermayer (1999) implies the existence of a process $Y \in \mathcal{Y}(1)$ [the set $\mathcal{Y}(y)$ has been defined in (10)] such that $Y_T \geq h$. Taking into account (27), we deduce that the process $Z$ defined by

$$Z_t = \begin{cases} Y_t, & t < T, \\ h, & t = T, \end{cases}$$

belongs to $\mathcal{Y}(1,p)$. Therefore, $h \in \mathcal{D}(1,p)$. □

For the proof of Theorem 1 we shall also need the following lemma.

LEMMA 11. *Let $\mathcal{E}$ be a set of nonnegative random variables which is convex and contains a strictly positive constant. Then*

(33) $$\sup_{g \in \mathcal{E}} \mathbb{E}[U(xg)] = \sup_{g \in \operatorname{cl}\mathcal{E}} \mathbb{E}[U(xg)], \qquad x > 0,$$

*where $\operatorname{cl}\mathcal{E}$ denotes the closure of $\mathcal{E}$ with respect to convergence in probability.*

PROOF. Without loss of generality, we can assume that the set $\mathcal{E}$ contains 1. Denote, for $x > 0$,

$$\phi(x) \triangleq \sup_{g \in \mathcal{E}} \mathbb{E}[U(xg)], \qquad \psi(x) \triangleq \sup_{g \in \operatorname{cl}\mathcal{E}} \mathbb{E}[U(xg)].$$

Clearly, $\phi$ and $\psi$ are concave functions and $\phi \leq \psi$. If $\phi(x) = \infty$ for some $x > 0$, then, due to concavity, $\phi$ is infinite for all arguments and the assertion of the lemma is trivial. Hereafter we assume that $\phi$ is finite.

Fix $x > 0$ and $g \in \operatorname{cl}\mathcal{E}$. Let $(g^n)_{n \geq 1}$ be a sequence in $\mathcal{E}$ that converges to $g$ almost surely. For any $\delta > 0$, we have

$$\mathbb{E}[U(xg)] \leq \mathbb{E}[U(xg + \delta)] \leq \liminf_{n \to \infty} \mathbb{E}[U(xg^n + \delta)] \leq \phi(x + \delta),$$

where the first inequality holds true because $U$ is increasing, the second one follows from Fatou's lemma and the third one follows from the facts that $\mathcal{E}$ is convex and contains 1. Since $\phi$ is concave, it is continuous. It follows that

$$\psi(x) = \sup_{g \in \operatorname{cl}\mathcal{E}} \mathbb{E}[U(xg)] \leq \lim_{\delta \to 0} \phi(x + \delta) = \phi(x). \qquad \square$$

PROOF OF THEOREM 1. The function $u$ is clearly concave on $\mathcal{K}$. Since the set $\mathcal{K}$ is open and $u(x,q) < \infty$ for some $(x,q) \in \mathcal{K}$, we deduce that $u$ is finitely valued on $\mathcal{K}$.

For $(y,r) \in \mathcal{L}$, we define the sets

$$A(y,r) \triangleq \{(x,q) \in \mathcal{K} : xy + \langle q, r \rangle \leq 1\},$$
$$\widetilde{\mathcal{C}} \triangleq \bigcup_{(x,q) \in A(y,r)} \mathcal{C}(x,q),$$

and denote by $\mathcal{C}$ the closure of $\widetilde{\mathcal{C}}$ with respect to convergence in probability. From Lemma 11 and the definition of the set $\widetilde{\mathcal{C}}$, we deduce that

$$\sup_{g \in \mathcal{C}} \mathbb{E}[U(zg)] = \sup_{g \in \widetilde{\mathcal{C}}} \mathbb{E}[U(zg)] = \sup_{(x,q) \in zA(y,r)} u(x,q), \qquad z > 0.$$

Further, we claim that

(34) $$\sup_{(x,q) \in zA(y,r)} u(x,q) < \infty, \qquad z > 0.$$

To prove (34), assume that the set $\mathcal{L}$ is open. As was explained in Remark 6, this is without any loss of generality. In this case, the set $A(y,r)$ is bounded and (34) follows from the concavity of $u$.

From Proposition 1, we deduce that

$$h \in \mathcal{D}(y,r) \quad \Leftrightarrow \quad \mathbb{E}[gh] \leq 1 \quad \forall g \in \mathcal{C}.$$

It follows that the sets $\mathcal{C}$, and $\mathcal{D}(y,r)$ satisfy the assumptions of Theorem 3.1 in Kramkov and Schachermayer (1999). If $v(y,r) < \infty$, then this theorem implies the existence and uniqueness of the solution $\widehat{h}(y,r)$ to (25). It also implies the existence of a strictly positive constant $c = c(y,r)$ such that $v(cy, cr) < \infty$ as well as the second conjugacy relationship in (14):

(35) $$v(y,r) = \sup_{z>0} \left\{ \sup_{g \in \mathcal{C}} \mathbb{E}[U(zg)] - z \right\} = \sup_{(x,y) \in \mathcal{K}} \{u(x,q) - xy - \langle q, r \rangle\}.$$

Finally, the first equation in (14) follows from the second one by the general properties of conjugate functions [see, e.g., Rockafellar (1970), Section 12] and the fact that the relative interior of the domain of the function $v$ defined by (35) is a subset of $\mathcal{L}$. □

For the proof of Theorem 2, we shall also need the following lemma.

LEMMA 12. *Assume that the assumptions of Proposition* 1 *hold true. Let sequences* $(y^n, r^n) \in \mathcal{L}$ *and* $h^n \in \mathcal{D}(y^n, r^n)$, $n \geq 1$, *converge to* $(y,r) \in \mathbf{R}^{N+1}$ *and* $h \in \mathbf{L}_+^0$, *respectively. If* $h$ *is a strictly positive random variable, then* $(y,r) \in \mathcal{L}$ *and* $h \in \mathcal{D}(y,r)$.

PROOF. Without loss of generality, we can assume that the assertions of Lemma 7 are valid. Let $(x, q) \in \operatorname{cl} \mathcal{K}$ and $X \in \mathcal{X}(x, q)$ be as in assertion (iii) of Lemma 7, that is, such that

$$\mathbb{P}[X_T + \langle q, f \rangle > 0] > 0.$$

Combining this property with Proposition 1 and Fatou's lemma, we deduce that

$$0 < \mathbb{E}[h(X_T + \langle q, f \rangle)] \leq xy + \langle q, r \rangle.$$

As $(x, q)$ is an arbitrary element of $\operatorname{cl} \mathcal{K}$, this implies that $(y, r) \in \mathcal{L}$. Finally, Fatou's lemma and Proposition 1 imply that $h \in \mathcal{D}(y, r)$. □

PROOF OF THEOREM 2. Since

$$(36) \qquad U(x) \leq V(y) + xy, \qquad x > 0, y > 0,$$

we deduce from Proposition 1 that if $v(y, r) < \infty$ for some $(y, r) \in \mathcal{L}$, then $u(x, q) < \infty$ for all $(x, q) \in \mathcal{K}$. In particular, all of the assumptions of Theorem 1 hold true. For $(x, q) \in \mathcal{K}$, define the sets

$$B(x, q) \triangleq \{(y, r) \in \mathcal{L} : xy + \langle q, r \rangle \leq 1\},$$

$$\widetilde{\mathcal{D}} \triangleq \bigcup_{(y,r) \in B(x,q)} \mathcal{D}(y, r)$$

and denote by $\mathcal{D}$ the closure of $\widetilde{\mathcal{D}}$ with respect to convergence in probability. Using these definitions, we clearly have

$$\inf_{h \in \mathcal{D}} \mathbb{E}[V(zh)] \leq \inf_{h \in \widetilde{\mathcal{D}}} \mathbb{E}[V(zh)] = \inf_{(y,r) \in zB(x,q)} v(y, r) < \infty, \qquad z > 0.$$

From Proposition 1, we deduce that

$$g \in \mathcal{C}(x, q) \quad \Leftrightarrow \quad \mathbb{E}[gh] \leq 1 \qquad \forall h \in \mathcal{D}.$$

It follows that the sets $\mathcal{C}(x, q)$ and $\mathcal{D}$ satisfy the assumptions of Theorem 4 in Kramkov and Schachermayer (2003). This implies the existence and uniqueness of the solutions $\widehat{g}(x, q)$ to (24) and $\widehat{X}(x, q)$ to (7), which are related by the following equality:

$$\widehat{X}_T(x, q) = \widehat{g}(x, q) - \langle q, f \rangle.$$

Denoting

$$\widehat{h} \triangleq U'(\widehat{g}(x, q)), \qquad z \triangleq \mathbb{E}[\widehat{g}(x, q) \widehat{h}],$$

we also deduce from this theorem that $\widehat{h}$ belongs to $z\mathcal{D}$ and is the unique solution to the optimization problem:

$$\mathbb{E}[V(\widehat{h})] = \inf_{h \in z\mathcal{D}} \mathbb{E}[V(h)].$$

Since $\mathbb{P}[\widehat{h} > 0] = 1$, we deduce from Lemma 12 the existence of $(y,r) \in B(x,q)$ such that the set $\mathcal{D}(y,r)$ contains $\widehat{h}$. Clearly, for this $(y,r)$ we have that

$$xy + \langle q, r \rangle = z = \mathbb{E}[\widehat{g}(x,q)\widehat{h}\,] \tag{37}$$

and that the solution to (25) satisfies

$$\widehat{h}(y,r) = \widehat{h} = U'(\widehat{g}(x,q)). \tag{38}$$

Since

$$U(\widehat{g}(x,q)) = V(\widehat{h}) + \widehat{g}(x,q)\widehat{h},$$

we deduce that

$$u(x,q) = v(y,r) + xy + \langle q, r \rangle, \tag{39}$$

which, according to Rockafellar (1970), Theorem 23.5, is equivalent to $(y,r) \in \partial u(x,q)$. In particular, we have

$$\partial u(x,q) \cap \mathcal{L} \neq \varnothing. \tag{40}$$

Conversely, if $(y,r)$ belongs to $\mathcal{L}$ and satisfies (39), that is, $(y,r) \in \mathcal{L} \cap \partial u(x,q)$, then

$$\mathbb{E}[|V(\widehat{h}(y,r)) + \widehat{g}(x,q)\widehat{h}(y,r) - U(\widehat{g}(x,q))|]$$
$$= \mathbb{E}[V(\widehat{h}(y,r)) + \widehat{g}(x,q)\widehat{h}(y,r) - U(\widehat{g}(x,q))]$$
$$\leq v(y,r) + xy + \langle q, r \rangle - u(x,q) = 0,$$

which readily implies (37) and (38).

To conclude the proof of the theorem we now need to show that

$$\partial u(x,q) \subset \mathcal{L}.$$

Let $(y,r) \in \partial u(x,q)$. From (40) and the fact that $\partial u(x,q)$ is a closed convex set, we deduce the existence of a sequence $(y_n, r_n)$ in $\partial u(x,q) \cap \mathcal{L}$ that converges to $(y,r)$. As we have proved already, any of the sets $\mathcal{D}(y_n, r_n)$ contains the strictly positive random variable $U'(\widehat{g}(x,q))$. Lemma 12 implies now that $(y,r) \in \mathcal{L}$. □

**Acknowledgments.** We thank Pierre Collin Dufresne, Nizar Touzi, Mihai Sirbu and Walter Schachermayer for discussions on the topic of this paper.

HEC Montréal and CIRANO  
3000 côte Ste Catherine  
Montréal, Quebec H3T 2A7  
Canada  
e-mail: julien.hugonnier@hec.ca

Department of Mathematical Sciences  
Carnegie Mellon University  
5000 Forbes Avenue  
Pittsburgh, Pennsylvania 15213  
USA  
e-mail: kramkov@andrew.cmu.edu